\documentclass[12pt]{amsart}
\usepackage{bm,hyperref,fullpage,color}
\usepackage[showonlyrefs]{mathtools}

\def\nablaIT{\setbox0=\hbox{$\nabla$}%
	\pdfliteral{q 1 0 .3 1 0 0 cm}\rlap{$\nabla$}\pdfliteral{Q}\kern\wd0 }
\newcommand\itnabla{\operatorname{\nablaIT}{}}

\title{An interpretation of \\ Temam's extra force in  the \\ quasi-incompressible Navier-Stokes system}
\author{Giuseppe Tomassetti${}^*$}
\thanks{${}^*$Universit\`a degli Studi Roma Tre, Dipartimento di Ingegneria, Via Vito Volterra 62, 00154 Roma, Italy. Email: \texttt{giuseppe.tomassetti@uniroma3.it}.}

\begin{document}
\begin{abstract}
 We discuss the role of the  extra force density $$\mathbf f_{\rm e}=-\frac 1 2(\nabla\cdot\mathbf v)\mathbf v$$ in the adimensionalized system of partial differential equations
  \begin{equation*}
    \left\{
    \begin{aligned}
      &\frac{\partial\mathbf v}{\partial \rm t}+(\mathbf v\cdot\nabla)\mathbf v+\nabla \mathrm p-\frac 1 {\rm Re}\mathit{\Delta}\mathbf v=\mathbf f+\mathbf f_{\rm e},\\
      &\frac 1 {\mathrm K}\frac{\partial \mathrm p}{\partial \rm t}+\nabla\cdot\mathbf v=0,\qquad \mathrm K>>1,
      \end{aligned}\right.
  \end{equation*}
  whose weak solution, with appropriate initial and boundary conditions, has been proved in \href{https://doi.org/10.1007/BF00247678}{[Arch. Rat. Mech. Analysis, 32:135--153]} to preserve the balance of energy while approximating, in the limit $\mathrm K\to\infty$, the weak solution of the incompressible Navier-Stokes system, where the extra force vanishes. Taking the cue from \href{https://doi.org/10.1007/BF01782609}{[Ann. Mat. Pura Appl. 172:103--124]}, we provide a mechanical interpretation of the extra force density $\mathbf f_{\rm e}$, arguing that it is a manifestation of inertia. 
\end{abstract}

 \maketitle


\section{Introduction}
The quasi-incompressible Navier-Stokes (N--S) system is an approximation of the fully compressible N--S system which may be used in those circumstances when the bulk elastic modulus is very large, in comparison with the occuring pressures, but not infinite; it represents a trade-off between the fully compressible model, which is known to be very complicated from the analytical and numerical point of view, and the incompressible model, which cannot describe several interesting phenomena, such as for instance the propagation of pressure waves \cite{roubicek2020}. In addition, the quasi-incompressible version of the N--S system makes it possible to obtain an a priori control of the pressure, which is essential both in the analysis and the numerics of fluids with pressure-dependent viscosity \cite{hron,hron2}. Further applications, more recently, are concerned with models for geodynamical flows \cite{roubicektomassetti2020}. Such wide spectrum of applications motivates a deeper scrutiny of the quasi-incompressible N-S system.
\smallskip

The quasi-incompressible N--S system was originally introduced by R. Temam \cite{temam1968} (see also the article \cite{Temam} and the monography \cite{Temam}) as a mathematical device for the approximation of the incompressible N--S system, to alleviate the stiffness of the incompressibility constraint \cite{prohl1997}. In Temam's formulation, an extra force appears, whose role from the purely mathematical standpoint is to guarantee the validity of an estimate which corresponds to the balance of energy. Although this idea has generated a substantial amount of mathematical work in recent years \cite{bersellispirito2018,donatelli2010,donatellimarcati2006,donatellispirito2011,oskolkov1973,roubicek2020,shao2017,zhao2013}, from the engineering point of view it would be interesting to find some interpretation of  Temam's extra force.
\smallskip
 
In this note we  propose a mechanical interpretation of this extra force, based on a cue from \cite{Podio97}. Our line of argument is the following: 1) by invoking d'Alembert's Principle, we write the balance of momentum as an equilibrium equation involving inertial and non-inertial forces; 2) we regard the specification of inertial forces as a \emph{postulate} subjected to the fundamental requirement that the \emph{inertial power}, {\emph i.e.}, the power expended by the inertial forces, be equal to the negative rate of change of kinetic energy; 3) we argue that if mass density is assumed \emph{constant}, but flow has \emph{non-null divergence}, then an additional term must be added to the standard expression of the inertial force to maintain  consistency with the requirement in 2). In our opinion, this approach may open the way to several generalizations. For example, an open problem is the modeling of layered or stratified continua (such as in geophysics) where density varies with depth. We also mention problems involving growing bodies \cite{lubardahoger} or bodies with self diffusion \cite{epsteingoriely}. In these problems conservation of mass does not hold in general, and choosing the appropriate form of the inertial force is a delicate matter.
\smallskip

This note is organized as follows. Section \ref{sec:compr-navi-stok} contains a brief account concerning the derivation of the Navier-Stokes equation. Its main purpose is to establish a notational framework and to illustrate our point of view concerning the role of inertial forces, which, in the spirit of Noll \cite{noll1963}, we regard as the object of a prescription. Section \ref{sec:incompr-navi-stok2} is devoted to the incompressible NS system, whereas Section \ref{sec:quasi-incompr-navi} is concerned with the quasi-incompressible approximation proposed in \cite{Temam}. The most important part of this note is Section \ref{sec:mech-interpr-stab} where, adopting the Lagrangian point of view, we show that the extra term introduced in \cite{Temam} guarantees the balance between the time derivative of the specific --- \emph{i.e.} per unit reference volume --- kinetic energy and the specific power expended by inertial forces. We conclude with some remarks in Section \ref{sec:conclusions}.
 
\section{The compressible Navier-Stokes system}\label{sec:compr-navi-stok}
The Navier--Stokes (N--S) equation, whose most popular version is
\begin{equation}\label{eq:1}\tag{$\rm NS$}
  \varrho\frac{\partial\bm v}{\partial t}+\varrho(\bm v\cdot\itnabla)\bm v+\itnabla p-\mu\mathit{\Delta}\bm v-\Big(\zeta+\frac \mu 3\Big)\itnabla(\itnabla\cdot\bm v)=\bm f,
\end{equation}
constitutes a basic model in \emph{fluid dynamics}. In this expression, the spatio-temporal fields $\varrho(x,t)$ and $\bm v(x,t)$ are, respectively, the \emph{mass density} $({\rm kg}/{\rm m}^{3})$ and the \emph{velocity} $({\rm m}/{\rm s})$; $p(x,t)$ is the \emph{pressure} field $({\rm Pa})$, and $\bm f(x,t)$ is the \emph{non-inertial body force} field $(\rm N/m^3)$. 
Here we use the format typical in fluid mechanics. In particular, the term $(\bm v\cdot\itnabla)\bm v=(\itnabla\boldsymbol v)\boldsymbol v$, that is, the application of the linear operator $\itnabla\boldsymbol v$ to $\boldsymbol v$. The terms $\mathit{\Delta}\bm v$ and $\itnabla(\itnabla\cdot\bm v)$ are, respectively, the Laplacian of the velocity field and the gradient of its divergence. These terms carry dissipative effects proportional to the \emph{dynamic viscosity} $\mu$ $(\rm Pa/s)$ and the \emph{volumetric viscosity} $\zeta$, assumed to be constant. 


A possible way to obtain the N-S system, in the spirit of D'Alembert's principle, is to start from the pointwise equilibrium equation
\begin{equation}\label{eq:21}
  \bm f_{\rm i}+\bm f+\itnabla\cdot\boldsymbol\sigma=\bm 0,
\end{equation}
which dictates that the inertial body force $\bm f_{\rm i}$, the (prescribed) non-inertial body force $\bm f$, and the system of internal forces be in equilibrium, the  system of internal forces being equipollent to a body force whose density is the divergence of the Cauchy stress $\boldsymbol\sigma$. In  particular, \eqref{eq:1} follows from the prescription (see \cite[Eq. 19.48]{gurtinfriedanand})
\begin{equation}\label{eq:11}
 \bm f_{\rm i}=-\varrho\Big(\frac{\partial\bm v}{\partial t}+(\bm v\cdot\itnabla)\bm v\Big),
\end{equation}
and
\begin{equation*}
  \boldsymbol\sigma=-p\bm I+\mu\bigg(\itnabla\bm v+\itnabla\bm v^T+\Big(\frac {\zeta}\mu-\frac 2 3\Big) (\itnabla\cdot\bm v)\bm I\bigg)
\end{equation*}
for the inertial force and for the Cauchy stress, respectively. The definition \eqref{eq:11} guarantees (see for example \cite[Chap. 19]{gurtinfriedanand}) that, given any spatial region convecting with the body, the rate of change of the total linear momentum in that region be equal to minus the integral of the inertial body force $\bm f_{\rm i}$ over that region. In this argument, balance of mass plays a key role. On the other hand, several authors, starting with Noll \cite{noll1963} (see also \cite[p. 144]{gurtinfriedanand} and \cite[Sec. 2]{Podio97}) consider \eqref{eq:11} as a postulate.

The mass-density and velocity fields are not independent, being related by the \emph{mass-balance equation}:
\begin{equation}\label{eq:2}\tag{{\rm MBE}}
  \frac{\partial\varrho}{\partial t}+\itnabla\cdot(\varrho\bm v)=0.
\end{equation}
The system \eqref{eq:1}-\eqref{eq:2} is  closed by the \emph{Equation of State}
\begin{equation}\label{eq:4}\tag{\rm ES}
  p=\widehat p(\varrho),
\end{equation}
which makes explicit the dependence of \emph{thermodynamic pressure} $p$ on mass density through the \emph{constituive response function} $\widehat p$. We refer to the system \eqref{eq:1}, \eqref{eq:2}, \eqref{eq:4} as the \emph{compressible Navier--Stokes system}. As pointed out in \cite{rajagopal2015}, the use of the word pressure is not free from ambiguities. In what follows, we shall use the noun ``pressure'' as a reference to the thermodynamic pressure even if the true mechanical pressure is indeed $p_{\rm m}=-\frac 13 \operatorname{tr}\boldsymbol\sigma$. The two pressures coincide when $\zeta=0$, an assumption known as Stokes' hypothesis. This hypothesis is known to hold only for monoatomic gases \cite{tisza1942}, and whose untenability has been recently discussed in \cite{rajagopal2013}.

\section{The incompressible Navier-Stokes system and its approximation}\label{sec:incompr-navi-stok2}
As mentioned in the introduction, for fluids having high bulk modulus, the density can be considered  equal to a constant $\varrho_*$, which may be taken, for example, to be the density measured at the atmospherical pressure $p_*$. If we replace the Equation of State \eqref{eq:4} with the constraint
\begin{equation}
  \varrho(x,t)=\varrho_*,
\end{equation}
then pointwise equilibrium and mass balance yield the \emph{incompressible Navier-Stokes system}:
\begin{equation}\label{eq:5}
  \left\{
  \begin{aligned}
  &\varrho_*\frac{\partial\bm v}{\partial t}+\varrho_*(\bm v\cdot\itnabla)\bm v+\itnabla p-\mu\mathit{\Delta}\bm v=\bm f,\\
  &\itnabla\cdot\bm v= 0,
\end{aligned}
\right.
\end{equation}
with $p$ being now a \emph{reactive pressure}, {\rm i.e.}, not anymore determined by an equation of state, yet needed to enforce the incompressibility constraint $\itnabla\cdot\bm v=0$. 

When computing numerical solutions of the incompressible N--S system \eqref{eq:5}, the incompressibility constraint poses several issues \cite{prohl1997}. One possible approach to alleviate these issues is to \emph{relax the constraint}. A possible relaxation scheme involves using the compressible N--S system with a reasonably simple choice for the constitutive response function $\widehat p$:
\begin{equation}\label{eq:3}
  \widehat p(\varrho)=K\Big(\frac \varrho {\varrho_*}-1\Big)+p_*,
\end{equation}
where $K$, the \emph{bulk modulus} $[{\rm Pa}]$, is choosen as large as it takes to nullify, within the required numerical precision, the departure of $\varrho$ from its reference value $\varrho_*$. When \eqref{eq:3} is adopted, the inversion of \eqref{eq:4} yields
\begin{equation}\label{eq:6}
  \varrho=\varrho_*\Big(1+\frac {p-p_*} {K}\Big)=:\widehat\varrho(p/K).
\end{equation}
On substituting \eqref{eq:6} into \eqref{eq:1} and \eqref{eq:2} we obtain a system with respect to the unknowns $\bm v$ and $p$:
\begin{equation}\label{eq:10}
  \left\{
  \begin{aligned}
    &\widehat\varrho(p/K)\frac{\partial\bm v}{\partial t}+\widehat\varrho(p/K)(\bm v\cdot\itnabla)\bm v+\itnabla p-\mu\mathit{\Delta}\bm v-\Big(\zeta+\frac \mu 3\Big)\itnabla(\itnabla\cdot\bm v)=\bm f,\\
    &\frac{\varrho_*} K\frac{\partial p}{\partial t}+\itnabla\cdot(\widehat\varrho(p/K)\bm v)=0.
  \end{aligned}
\right.
\end{equation}
In view of \eqref{eq:6} one would expect that if the bulk modulus tends to infinity, and if pressure remains uniformly bounded in norm by a constant, then mass density should tend to the reference value $\varrho_*$:
\begin{equation*}
  K\to\infty \quad\text{and}\quad|p-p_*c|<\text{const.}\qquad \Rightarrow \qquad\varrho\to\varrho_*,
\end{equation*}
so much so the incompressibility constraint would be recovered:
\begin{equation*}
  \itnabla\cdot\bm v\to 0.
\end{equation*}
For the rest of this discussion, we shall find it useful to write the incompressible Navier-Stokes system \eqref{eq:5} in the dimensionless form:
\begin{equation}\label{eq:5bis}
  \left\{
  \begin{aligned}
  &\frac{\partial\mathbf v}{\partial \rm t}+(\mathbf v\cdot\nabla)\mathbf v+\nabla{\rm p}-\frac{1}{\rm Re}\Delta\mathbf v=\mathbf f,\\
  &\nabla\cdot\mathbf v= 0,
\end{aligned}
\right.
\end{equation}
where  ${\rm Re}=\varrho_* LV/\mu$ is the Reynolds Number and upright fonts denote dimensionless fields, defined by
\begin{equation}\label{eq:14}
  \begin{aligned}
    \mathbf v(\mathbf x,\mathrm t)=V\boldsymbol v(x,t),\qquad \mathrm{p}(\mathbf x,\mathrm t)=\frac{p(x,t)-p_*}{\varrho_*V^2},\qquad \mathbf f=\frac{L}{\varrho_* V^2}\bm f,\qquad \mathbf x=\frac {\boldsymbol x}L,\qquad \mathrm t=\frac V L  t,
    \end{aligned}
\end{equation}
with  $V$ and $L$ a characteristic velocity and a characteristic length scale, and the symbols $\nabla$ and $\Delta$ denote, respectively, the gradient and the Laplacian with respect to the dimensionless space variables. 

\section{The quasi-incompressible Navier-Stokes system}\label{sec:quasi-incompr-navi}
The incompressible N-S system \eqref{eq:5} is simpler than the fully compressible system \eqref{eq:10}. However, it has the disadvantage that pressure is a reactive field, not determined constitutively. As a side effect, it leaves out of the picture several physical phenomena, such as for instance the propagation of pressure waves. This fact limits the range of applications of the incompressible model.

A trade-off between full compressibility and quasi-incompressibility is the notion of quasi-compressibility due to Temam. Using the dimensionless quantities introduced in the previous section, Temam' system reads: 
\begin{equation}\label{eq:9}
  \left\{
  \begin{aligned}
  &\frac{\partial\mathbf v}{\partial\mathrm t}+(\mathbf v\cdot\nabla)\mathbf v+\nabla {\rm p}-\frac{1}{\rm Re}\Delta\mathbf  v= \mathbf  f-\frac 1 2(\nabla\cdot\mathbf  v)\mathbf  v,\\
  &\frac{1}{\mathrm K}\frac{\partial\mathrm p}{\partial \mathrm t}+\nabla\cdot\mathbf  v=0,
  \end{aligned}\right.
\end{equation}
The heuristics behind \eqref{eq:9} is that in some intermediate regime when the bulk modulus is large, but not infinite, and pressure remains bounded, then by \eqref{eq:6} the density will tend to a constant $\varrho_{*}$, but still variations of the pressure may be captured by the second equation of \eqref{eq:10}. A derivation of the notion of quasi-incompressibility in the context of the theory of mixtures may be found in \cite{feireisl}, where a quasi-incompressible fluid is defined as a mixture consisting of two constituents, under the assumption of mass additivity and volume additivity constraints, and that the true densities of the constituents are constant. System \eqref{eq:9} is obtained from \eqref{eq:10} through the following \emph{four steps}:

(i) approximate the \emph{actual mass density} $\varrho=\widehat\varrho(p)$ with its reference value $\varrho_*$, i.e., perform the formal substitution:
\begin{equation*}
  \widehat\varrho(p/K)\mapsto\varrho_*;
\end{equation*}

(ii) dispense of the volumetric term $-(\zeta+\mu/3)\nabla(\nabla\cdot\bm v)$;

(iii) pass to dimensionless variables by performing the substitutions in \eqref{eq:14}, and by introducing the \emph{dimensionless bulk modulus}
\[\mathrm K=\frac{K}{\varrho_* V^2};
\]

(iv) add the adimensionalized \emph{extra force density}
\begin{equation}\label{eq:12}
  \mathbf f_{\rm e}=-\frac 1 2 (\nabla\cdot\mathbf  v)\mathbf  v
\end{equation}
to the adimensionalized non-inertial body force $\mathbf f$.

One of the advantages of \eqref{eq:9} over \eqref{eq:10} is that the term that contains the time derivative of the velocity  is \emph{linear} with respect to the unknowns ($\mathbf v$ and $\mathrm p$), which is a desirable feature from the point of view of both mathematical and numerical analysis. Mathematically, the adoption of \eqref{eq:9} is justified by the convergence theorem proved in \cite{Temam}. We attempt in the next section to provide a mechanical reading of this extra term, based on the notion of kinetic energy.

\section{A mechanical interpretation of the  extra force}\label{sec:mech-interpr-stab}
We show in this section that, besides analytical convenience, there is a somehow deeper (from the standpoint of  mechanics) motivation  for the introduction of the extra force \eqref{eq:12} in the first of \eqref{eq:9}. For the reader's convenience we split our presentation in two parts. In the first part we recapitulate the procedure leading to the standard prescription \eqref{eq:11} for the inertial-force density. In the second part we offer our interpretation of the extra force. 
\subsection{The standard prescription for the inertial force}
First, we assume that the \emph{spatial density of  kinetic energy} be
\begin{equation}\label{eq:8}
  \kappa=\frac \varrho 2 |\bm v|^2.
\end{equation}
Next, we choose any reference configuration where the mass density is \emph{a constant and uniform} $\varrho_*$. Then, the \emph{Jacobian} of the deformation map which takes the reference configuration into the current configuration at point $x$ and time $t$ satisfies
\begin{equation}\label{eq:7}
  J(x,t)=\frac{\varrho_*}{\varrho(x,t)}.
\end{equation}
Therefore, the kinetic energy per unit \emph{referential} volume is
\begin{equation}\label{eq:19}
    \kappa_{\rm r}=J\kappa=\frac {J\varrho}2 |\bm v|^2=\frac {\varrho_*}2 |\bm v|^2.
  \end{equation}
Thus, the material time derivative of the referential kinetic energy is:
  \begin{equation}
    \dot\kappa_{\rm r}=\varrho_*\dot{\bm v}\cdot\bm v,
  \end{equation}
  where $\dot{\bm v}$ is the material time derivative of $\bm v$.

  Now, following \cite{Podio97}, we require that the rate of change of specific (\emph{i.e.} per unit referential volume or, equivalently, per unit mass) kinetic energy plus the specific inertial power be null during \emph{every possible motion}. In the present case, this requirement takes the form:
  \begin{equation}
    (\bm f_{\rm i,r}+\varrho_*\dot{\bm v})\cdot\bm v=0.
  \end{equation}
This requirement singles out the inertial force up to a powerless contribution which, when taken to be null, prompts the following choice for the referential inertial force:
  \begin{equation}\label{eq:20}
    \bm f_{\rm i,r}=-\varrho_*\dot{\bm v},
  \end{equation}
whence the following  prescription for the inertial force density in the current configuration:
  \begin{equation}\label{eq:16}
    \bm f_{\rm i}=J^{-1}\bm f_{\rm i,r}=\frac\varrho{\varrho_*}\bm f_{\rm i,r}=-\varrho\dot{\bm v}.
  \end{equation}
  The step from \eqref{eq:16} to \eqref{eq:11} is immediate on  recalling that
  \[
    \dot{\bm v}=\frac{\partial\bm v}{\partial t}+(\bm v\cdot\itnabla)\bm v.
    \]
\subsection{A non standard prescription for the inertial force}
Suppose that if instead of \eqref{eq:8} we make the following  choice for the \emph{spatial} (\emph{i.e.} in the current configuration) density of kinetic energy:
\begin{equation}\label{eq:13}
  \kappa^*=\frac {\varrho_*}2 |\bm v|^2.
\end{equation}
Then, as we show in the foregoing, the procedure outlined in the previous subsection leads to the following  prescription for the spatial inertial-force density:
\begin{equation}\label{eq:15}
  \bm f_{\rm i}^*=-\varrho_*\Bigl(\frac{\partial\bm v}{\partial t}+(\bm v\cdot\itnabla)\bm v+\frac {1} 2 (\itnabla\cdot\bm v)\bm v\Bigr).
\end{equation}
Based on this result we argue that, when performing Step (i) in the procedure that leads from \eqref{eq:10} to \eqref{eq:9}, the replacement of the actual mass density with its reference value $\varrho_*$ {\bf must be accompanied by the adoption of \eqref{eq:15} in place of \eqref{eq:11} as the prescription for the inertial force density}. However, replacing \eqref{eq:11} with \eqref{eq:15} brings about the following additional term on the right-hand side of the equation:
\begin{equation}\label{eq:18}
  \bm f_{\rm e}=-\frac {\varrho_*} 2 (\itnabla\cdot\bm v)\bm v.
\end{equation}
When passing to dimensionless variables, this terms trasforms into the extra force $\mathbf f_{\rm e}$ which is added in Step (iv). Accordingly, we should interpret the extra term added in Step (iv) as a correction to the prescription of the inertial force consistent with the assumption that the kinetic energy has the form \eqref{eq:13}.

\subsection{Proof of \eqref{eq:15}.}
As a start, we recall that the mass-balance equation \eqref{eq:2} can be written as
\begin{equation}\label{eq:17}
  \dot\varrho+\varrho\itnabla\cdot\bm v=0.
\end{equation}
Next, on adopting \eqref{eq:13}, we obtain that the kinetic energy per unit \emph{referential} volume is
\begin{equation*}
    \kappa_{\rm r}^*=J\kappa^*=\frac {J\varrho_*}{2} |\bm v|^2\stackrel{\eqref{eq:7}}=\frac 12 \frac {\varrho_*^2}{\varrho} |\bm v|^2.
  \end{equation*}
Thus, the rate of change of the referential kinetic-energy density is:
  \begin{equation*}
    \dot\kappa_{\rm r}^*=\frac {\varrho_*^2}{\varrho}\dot{\bm v}\cdot\bm v-\frac {\dot\varrho}2 \frac {\varrho_*^2}{\varrho^2}|\bm v|^2\stackrel{\eqref{eq:17}}=\frac {\varrho_*^2}{\varrho}\dot{\bm v}\cdot\bm v+\frac {\varrho(\itnabla\cdot\bm v)}2 \frac {\varrho_*^2}{\varrho^2}|\bm v|^2=\frac {\varrho_*}{\varrho}\left(\varrho_*\dot{\bm v}+\frac {\varrho_*}2(\itnabla\cdot\bm v)\bm v\right)\cdot\bm v.
  \end{equation*}
 This leads us to the following choice for the referential density of the inertial force:
  \begin{equation*}
    \bm f_{\rm r}^*=-\frac {\varrho_*}{\varrho}\left(\varrho_*\dot{\bm v}+\frac {\varrho_*}2(\itnabla\cdot\bm v)\bm v\right),
  \end{equation*}
which yields, for the inertial force per unit volume in the current configuration, the expression
  \begin{equation*}
    \bm f_{\rm i}^*=\frac {\varrho}{\varrho_*}\bm f_{\rm i,r}^*=-\varrho_*\dot{\bm v}-\frac {\varrho_*}2(\itnabla\cdot\bm v)\bm v,
  \end{equation*}
  which coincides with \eqref{eq:15}, in view of \eqref{eq:12}.
\medskip

  In conclusion, we have shown that if the approximate expression \eqref{eq:13} is adopted for the spatial density of kinetic energy, where $\varrho_*$ is the mass density appearing in the incompressible Navier-Stokes system \eqref{eq:5}, then the extra force \eqref{eq:12} in the quasi-incompressible Navier-Stokes system \eqref{eq:9} emerges as a natural consequence of the requirement that the specific power expenditure of the inertial force be equal to minus the rate of change of the specific kinetic energy. In other words, the extra force  may be interpreted as  a  manifestation of inertia. 
 
  \section{Concluding remarks}\label{sec:conclusions} 
\paragraph{\bf Remark 1.} The format of the quasi-incompressible NS system, as studied by \cite{Temam}, involves a convective time derivative $\dot{\mathbf v}=\frac{\partial\mathbf v}{\partial t}+(\mathbf v\cdot\itnabla)\mathbf v$ in the first equation and the  partial derivative with respect to time in the second equation. If we had started from the following form of the mass-balance equation \eqref{eq:2}:
    \begin{equation*}
\dot{\varrho}+\varrho\itnabla\cdot\mathbf v=0,
\end{equation*}
then the application of our argument would have led us to:
\begin{equation*}
  \frac{1}{\mathrm K}\dot {\mathrm p}+\nabla\cdot\mathbf v=0.
\end{equation*}
This is indeed the equation used in \cite[Eq. (4b)]{roubicek2020}.

\medskip

  \paragraph{\bf Remark 2.} The standard prescription of the inertial force is invariant under a Galilean change of observer:
  \begin{equation}
    x^*=x+\bm wt,\qquad t^*=t+\tau,
  \end{equation}
  where $\bm w$ is a velocity and $\tau$ is a time. In fact, for $\bm v^*(x^*,t^*)=\bm v(x,t)+\bm w$, we have
  \begin{equation}
    \frac{\partial\bm v^*}{\partial t^*}=\frac{\partial\bm v}{\partial t}+(\bm w\cdot\itnabla)\bm v,\qquad \itnabla^*\bm v^*=\itnabla\bm v.
  \end{equation}
  Hence  
  \begin{equation}
    \frac{\partial\bm v}{\partial t}+(\bm v\cdot\itnabla)\bm v=\frac{\partial\bm v^*}{\partial t^*}-(\bm w\cdot\itnabla)\bm v+(\bm v^*\cdot\itnabla^*)\bm v^*+(\bm w\cdot\itnabla)\bm v=\frac{\partial\bm v^*}{\partial t^*}+(\bm v^*\cdot\itnabla^*)\bm v^*.
  \end{equation}
This is not the case with Temam's force. This observation uncovers a major drawback in the introuction of Temam's force when using the quasi-incompressible system to model slightly compressible fluids  and  opens up to the question concerning other possible types of regularization in the spirit of Temam's. Of course, different expressions for the extra force would arise upon  different choices of the kinetic energy. For instance, by adopting the standard choice of the kinetic-energy density \eqref{eq:8}, and by selecting $\varrho=\widehat\varrho(p/K)$ as given by the inversion \eqref{eq:6} of equation of state \eqref{eq:5}, we would recover the expression \eqref{eq:19} for the referential energy density, and hence we would get back to the first of \eqref{eq:10}, which then could be written as
  \begin{equation}
    \varrho_*\frac{\partial\bm v}{\partial t}+\varrho_*(\bm v\cdot\itnabla)\bm v-\itnabla p-\mu\Delta\bm v-\frac \mu 3\itnabla(\itnabla\cdot\bm v)=\bm f+\widetilde{\bm f}_{\rm e},
  \end{equation}
  with $\widetilde{\bm f}_{\rm e}=-\varrho_*\frac p K\big(\frac{\partial\bm v}{\partial t}+(\bm v\cdot\itnabla)\bm v\big)$ a properly Galiean-invariant force, which vanishes as $K$ tends to infinity if the absolute value of the pressure stays bounded.
  
In this respect, we notice that the quasi-incompressible system bears some resemblance with the Oberbeck--Boussinesq approximation for fluids that are mechanically incompressible but thermally compressible.  This system has been given a rigorous justification in \cite{rajagopalvergorisaccomandi2015} through asymptotic analysis. This prompts the question whether the quasi--incompressible system or a variant thereof may be obtained through a similar procedure.
\medskip

\paragraph{\bf Remark 3.} An alternative way to carry out the calculation in Section 5.3 is as follows. Let $J$ be any spatial field such that $\dot J=J\operatorname{div}\bm v$ (for example, $J$ may be the Jacobian of the deformation from some reference configuration). Then for every spatial field $\phi$ and for every  spatial region $\Omega_t$ convecting with the fluid (in the sense of \cite[p.~63]{gurtinfriedanand}), by Reynolds' transport theorem (\cite[p.~113]{gurtinfriedanand}), we have ${\rm d}/{\rm d}t\int_{\Omega_t}\phi$ = $\int_{\Omega_t}\big(\dot\phi+\phi\itnabla\cdot\bm v\big)$. In particular, taking $\phi$ =$ J^{-1}\varphi$, where $\varphi$ is another spatial field, we can write ${\rm d}/{\rm d}t\int_{\Omega_t}J^{-1}\varphi$ = $\int_{\Omega_t}\big(\dot{\overline{J^{-1}\varphi}}+J^{-1}\varphi\itnabla\cdot\bm v\big)$ = $\int_{\Omega_t}\big(-\frac 1 {J^2}\dot J\varphi+J^{-1}\dot\varphi+J^{-1}\varphi\itnabla\cdot\bm v\big)$ = $\int_{\Omega_t}J^{-1}\dot\varphi$. We apply this result with $\varphi=J\varrho_*|\bm v|^2/2$ to obtain ${\rm d}/{\rm d}t\int_{\Omega_t}\varrho_*|\bm v|^2/2$ = ${\rm d}/{\rm d}t\int_{\Omega_t} J^{-1}J\varrho_*|\bm v|^2/2$ = $\int_{\Omega_t}J^{-1}\dot{\overline{J\varrho_*|\bm v|^2/2}}$ = $\int_{\Omega_t}\varrho_*\bm v\cdot\dot{\bm v}+\varrho_*/2(\itnabla\cdot\bm v)\bm v\cdot\bm v$ = $\int_{\Omega_t}\big(\varrho_*{\dot{\bm v}}+\varrho_*/2(\itnabla\cdot\bm v)\bm v\big)\cdot\bm v$. Thus, one can see that if Temam's extra term $-\varrho_*/2(\itnabla\cdot\bm v)\bm v$ is added to the inertial force $-\varrho_*\dot{\bm v}$, then the power expended by the total inertial force is equal to minus the time derivative of the kinetic energy of any spatial region $\Omega_t$ convecting with the fluid.

\section{Acknowledgements}
The author is indebted to Alfredo Marzocchi, Paolo Podio-Guidugli, Pietro Prestininzi, and Rodolfo Repetto for their feedback. This work was supported by the Italian INdAM-GNFM (Istituto Nazionale di Alta Matematica -- Gruppo Nazionale per la Fisica Matematica) and the Grant of Excellence Departments, MIUR-Italy (Art.$1$, commi $314$-$337$, Legge $232$/$2016$). The author is also grateful to two anonymous reviewers, whose comments helped to rethink and improve the presentation.


\end{document}